\documentclass[11pt, lettersize, notitlepage, leqno, footskip=0.6cm]{article} 
\usepackage{amssymb}
\usepackage{amsthm}
\usepackage{amsmath}	
\usepackage{graphicx}
\usepackage{amscd}
	
\usepackage{indentfirst}

\usepackage{verbatim}
	
\usepackage{thmtools}
\usepackage[none]{hyphenat}	

\usepackage[bookmarks=true]{hyperref}
\usepackage{bookmark}

\newcommand{\bF}{\mathbb F}

\newcommand{\bz}{\mathbb Z}

\newcommand{\bq}{\mathbb Q}

\newcommand{\bFp}{\mathbb{F}_p}
\newcommand{\bFl}{\mathbb{F}_{\ell}}
\newcommand{\bFq}{\mathbb{F}_q}

\newcommand{\bFP}{\ov{\mathbb{F}}_p}

\newcommand{\bFQ}{\ov{\mathbb{F}}_q}

\newcommand{\vs}{\vspace{-2mm}}

\newcommand{\Gal}{\mathrm{Gal}}

\newcommand{\W}{\mathcal W}

\newcommand{\qbar}{\overline {\mathbb{Q}}}

\newcommand{\Mod}[1]{\ (\mathrm{mod}\ #1)}

\newcommand{\absq}{\mathrm{Gal}_{\bq}}

\newcommand{\lra}{\longrightarrow}

\newcommand{\hra}{\hookrightarrow}

\newcommand{\wti}{\widetilde}
\newcommand{\mf}{\mathfrak}
\newcommand{\mc}{\mathcal}

\newcommand{\mr}{\mathrm}

\newcommand{\mfp}{\mathfrak{p}}

\newcommand{\mfl}{\mathfrak{l}}

\newcommand{\zetamn}{\zeta_{mn}}

\newcommand{\al}{\alpha}
\newcommand{\be}{\beta}

\newcommand{\lamb}{\lambda}

\newcommand{\what}{\widehat}

\newcommand{\ov}{\overline}

\newcommand{\sub}{\subseteq}

\newcommand{\GalWB}{\Gal(\bq(\W_{B_{\pi}})/\bq)}

\DeclareMathOperator{\End}{End}

\DeclareMathOperator{\Mat}{Mat}

\DeclareMathOperator{\Nm}{Nm}
\DeclareMathOperator{\charac}{char}

\DeclareMathOperator{\IV}{IV}
\DeclareMathOperator{\Cl}{Cl}

\DeclareMathOperator{\isog}{isog}
\DeclareMathOperator{\Jac}{Jac}

\newtheorem{Thm}{Theorem}[section]
\newtheorem{Prop}[Thm]{Proposition}
\newtheorem{Lem}[Thm]{Lemma}
\newtheorem{Corr}[Thm]{Corollary}
\newtheorem{Algo}[Thm]{Algorithm}

\newtheorem{Def}{Definition}[section]

\declaretheoremstyle[%
  spaceabove=-2pt,%
  spacebelow=8pt,%
  headfont=\normalfont\itshape,%
  postheadspace=1em,%
  qed=\qedsymbol%
]{mystyle} 
\declaretheorem[name={Proof},style=mystyle,unnumbered,
]{prf}

\numberwithin{equation}{section}

\title{Abelian varieties with prescribed\\ embedding and full embedding degrees}
\author{Steve Thakur}  
\date{\vspace{-3ex}}

\begin{document} 
\maketitle

\begin{abstract} \noindent We study the problem of the embedding degree of an abelian variety over a finite field. In particular, we show that for a prescribed CM field $L$ of degree $\geq 4$, prescribed integers $m$, $n$ and a prescribed prime $\ell\equiv 1 \pmod{m\cdot n}$ that splits completely in $L$, there exists an ordinary abelian variety over a prime finite field with endomorphism algebra $L$, embedding degree $n$ with respect to $\ell$ and full embedding degree $m\cdot n$ with respect to $\ell$. We also study a class of absolutely simple higher dimensional abelian varieties whose endomorphism algebras are central over imaginary quadratic fields.\end{abstract}

\section{\fontsize{12}{12}\selectfont Introduction}

For a principally polarized abelian variety $A$ over a finite field $\bFq$ and a prime $\ell$ prime to $q$, the \textit{embedding degree} is the degree of the field extension $\bFq(\zeta_{\ell})$ obtained by attaching the $\ell$-th roots of unity. The \textit{full embedding degree} is the degree of the field extension $\bFq(A[\ell])$ obtained by attaching all $\ell$-torsion points of $A$.

The embedding and full embedding degrees for  varieties that are Jacobians of curves are natural notions in pairing-based cryptography. Such an abelian variety is principally polarized and is endowed with the non-degenerate \textit{Weil pairing} $$ A[\ell]\times A[\ell] \;\;\lra \;\; \mu_{\ell}  $$ from the subgroup scheme of $\ell$-torsion points to the $\ell$-th roots of unity in $\bFQ$. If $\bFq$ contains $\mu_{\ell}$, we also have the non-degenerate \textit{Tate pairing} $$ A[\ell](\bFq)\times A(\bFq)/\ell A(\bFq) \;\;\lra\;\; \bFq^*/(\bFq^*)^{\ell} .$$ 

In pairing-based cryptography, one typically chooses the prime $\ell$ and the embedding degree such that the discrete logarithm problems in the groups $A(\bFq)[\ell]$ (the $\ell$-torsion subgroup of the group of $\bFq$-valued points of $A$) and $\bFq(\mu_{\ell})^*$ are roughly of the same difficulty. This entails a large $\ell$ and a substantially smaller embedding degree. Abelian varieties with these properties are called \textit{pairing friendly}. Curves whose Jacobians allow for this are relatively rare and are also said to be pairing friendly curves.

The full embedding degree is an integer sandwiched in between the embedding degree and the prime $\ell$ and hence, one of the ways to look for suitable pairing-friendly abelian varieties is to look for abelian varieties with a prescribed embedding and full embedding degree. The crux of this short paper is some existence theorems related to this question. 

We only consider simple abelian varieties, i.e. those not isogenous to a product of two lower dimensional ableian varieties, since we can always reduce to this case. By Honda-Tate theory, the isogeny classes of simple abelian varieties over a finite field are in bijection with the Galois conjugacy classes of Weil numbers. Hence, the existence of abelian varieties with cetain properties is equivalent to the existence of suitable Weil numbers.

\subsection{\fontsize{11}{11}\selectfont Structure of the paper}

In section 1, we introduce some notations and background. We also state and prove a few lemmas that we will need repeatedly in the subsequent sections.

In section 2, we explore simple (albeit not \textit{absolutely} simple) supersingular abelian varieties as candidate for pairing friendly curves. 

In section 3, we discuss Koblitz's theorem on embedding degrees generalized to higher dimensional abelian varieties. We show that it fails to generalize to all but one class of absolutely simple abelian varieties of dimension two or higher.

In section 4, we study ordinary abelian varieties. In particular, we show that for a prescribed CM field $L$, prescribed integers $m$, $n$ and prescribed prime $\ell$ that splits completely in $L(\zetamn)$, there exists an ordinary abelian variety over a prime finite field with endomorphism algebra $L$ and with embeddng and full embedding degrees $n$ and $m\cdot n$ respectively with respect to the prime $\ell$.

In section 5, we study a class of absolutely simple abelian varieties with endomorphism algebra that are central over imaginary quadratic fields. This is the only class of higher dimensional abelian varieties where the notions of the embedding and the full embedding degrees coincide.

\subsection{\fontsize{11}{11}\selectfont Notations and background}

Throughout, $q$ denotes the size of a finite field which is the field of definition of an abelian variety under discussion. We denote the prime charactersitic of this field by $p$. Thus, $q$ is a $p$-power.

For an abelian variety $A$ over a field $F$, $\End^0(A)$ denotes the endomorphism algebra of $A$ over the field of definition. By Honda-Tate theory, we have the well-known bijection $$\{\text{Simple abelian varieties over } \bFq \text{ up to isogeny}\}\longleftrightarrow \{\text{Weil } q\text{-integers up to } \absq\text{-conjugacy} \}$$ induced by the map sending an abelian variety to its Frobenius. For a Weil number $\pi$, we write $B_{\pi}$ for the corresponding simple abelian variety (up to isogeny) over $\bFq$. The dimension of $B_{\pi}$ is given by $$2\dim B = [\bq(\pi):\bq]\cdot [\End^0(B_{\pi}):\bq(\pi)]^{1/2}.$$ Note that $\End^0(B_{\pi})$ is a central division algebra over $\bq(\pi)$ and hence, $[\End^0(B_{\pi}):\bq(\pi)]^{1/2}$ is an	integer. The characteristic polynomial of $B_{\pi}$ on the Tate representation $V_{\ell}(A):= T_{\ell}(A)\otimes_{\bz_{\ell}}\bq_l$ (for any prime ${\ell}\neq p$) is independent of ${\ell}$ and is given by \vspace{-1mm} $$P_{B_{\pi}}(X):= \prod\limits_{\sigma\in \absq} (X-\sigma(\pi))^{m_{\pi}}  $$ where $m_{\pi} = [\End^0(B_{\pi}):\bq(\pi)]^{1/2}.$ We denote by $\W_{B_{\pi}}$ the set of Galois conjugates of $\pi$. Hence, $\bq(\W_{B_{\pi}})$ is the splitting field of $P_{B_{\pi}}(X)$. The Galois group $\Gal(\bq(\W_{B_{\pi}})/\bq)$ is a subgroup of the wreath product $(\bz/2\bz)^g\rtimes S_g$. This is the same as the Galois group of the generic CM field of degree $2g$, where $g = [\bq(\pi+\ov{\pi}):\bq]$. 

The multiplicative group generated by $\W_{B_{\pi}}$ is denoted by $\Phi_{B_{\pi}}$. For a suitable base change $B_{\wti{\pi}}$ of $B_{\pi}$ to a larger extension, the group $\Phi_{B_{\wti{\pi}}}$ is a free abelian group. In the generic case, it is a free abelian group of rank $\dim B_{\pi}+1$.

For an abelian variety $B$ over a field $F$, $B(F)$ denotes the group of $F$-valued points on $B$. For an integer $n$, $B[n]$ denotes the $n$-torsion subgroup of $B(\ov{F})$ where $\ov{F}$ is the algebraic closure of $F$. The field extension generated by the coordinates of the points in $B[n]$ is denoted by $F(B[n])$.

For an integer $N\geq 1$, $\zeta_{_N}$ denotes a primitive $N$-th root of unity. $\phi(N)$ denotes the Euler totient function and $\Phi_N(X)$ denotes the $N$-th cyclotomic polynomial $$\Phi_N(X):= \prod\limits_{\substack{1\leq j\leq N \\ \gcd(j,N)= 1}} (X-\zeta_{_N}^j)\;\;\in\;\;\bz[X]$$ which is a degree $\phi(N)$ polynomial irreducible in $\bz[X]$. 

For a number field $K$, we denote its Hilbert class field by $\mr{H}(K)$. For a modulus $\mf{m}$ in a number field $K$, we denote the ray class field of $\mf{m}$ by $K^{\mf{m}}$.

We denote the Galois closure of $K$ over $\bq$ by $\wti{K}$. If $K$ is a CM field, we denote its maximal totally real subfield by $K^{+}$. We denote the ring of integers by $\mc{O}_K$.    

\begin{Def} \normalfont An abelian variety $A$ over a field $F$ is \textit{simple} if it does not contain a strict non-zero abelian subvariety. We say $A$ is \textit{absolutely} or \textit{geometrically} simple if the base change $A\times_F \ov{F}$ to the algebraic closure is simple.\end{Def}

\begin{Def} \normalfont An abelian variety $A$ is \textit{iso-simple} if it has a unique simple abelian subvariety up to isogeny.\end{Def}

\noindent We state a few well-known facts about abelian varieties over finite fields. We refer the reader to notes [Oo95] for proofs and further details.

\begin{Prop} For any simple abelian variety $B$ over a finite field $\bFq$, the abelian variety $B\times_{\bFq} \ov{\bF}_q$ is iso-simple. 
\end{Prop}
\vspace{-0.15cm}
With this setup, let $\pi$ be a Weil number corresponding to $B$ and let $\wti{B}$ be the unique simple component (up to isogeny) of the base change $B\times_{\bFq} \ov{\bF}_q$ to the algebraic closure. Let $N$ be the smallest integer such that $\wti{B}$ has a model over the field $\bF _{q^N}$. Then $\wti{B}$ corresponds to the Weil number $\pi^N$ and we have an isogeny \vspace{-0.15cm}$$B\times_{\bFq} \ov{\bF}_q =_{\isog} \wti{B}^{(\dim B)/N}.$$

\begin{Prop} Let $\pi$ be a Weil $q$-integer and write $D_{\pi}:= \End^0(B_{\pi})$. Then $D_{\pi}$ is a central division algebra over $\bq(\pi)$ and its Hasse invariants are given by $$\mr{inv}_{v}(D_{\pi}) = \begin{cases} 0 & \text{ if } v\nmid q.\\
\frac{1}{2} & \text{ if } v \text{ is real.}\\
[\bq(\pi)_{v}:\bq_p]\cdot \frac{v(\pi)}{v(q)} & \text{ if } v|q.\end{cases}$$\end{Prop}

\noindent In particular, $\End^0(B_{\pi})$ is commutative if and only if the local degrees $[\bq(\pi)_{v}:\bq_p]$ annihilate the Newton slopes $\frac{v(\pi)}{v(q)}$. For instance, if $B_{\pi}$ is ordinary, the slopes $\frac{v(\pi)}{v(q)}$ are either $0$ or $1$ and hence, $\End^0(B_{\pi})$ is commutative.

\subsection{\fontsize{11}{11}\selectfont A few preliminary results}

The following well-known result describes the embedding degree for a fixed prime $\ell$ and a simple abelian variety corresponding to a Weil number $\pi$ prime to $\ell$.

\begin{Prop} Let $B_{\pi}$ be a simple abelian variety over a finite field $\bFq$ corresponding to a Weil number $\pi$ and let $\ell$ be a prime not dividing $q$. The order $\big|B_{\pi}(\bFq)\big|$ of the group of $\bFq$-points is given by \vspace{-1mm} $$\big|B_{\pi}(\bFq)\big| = P_{B_{\pi}}(1) = \Nm(1-\pi)^{m_{\pi}}  \vspace{-1mm} $$ where $m_{\pi}$ is the multiplicity of $\pi$ in the characteristic polynomial. \end{Prop} 

The next few lemmas in this section might be obvious to most readers but we record them here since we will need them repeatedly in the subsequent sections.

\begin{Lem} Let $\pi$ be a Weil number and $B_{\pi}$ the corresponding abelian variety over a finite field $\bFq$. For any rational prime $\ell$, the following are equivalent:

\noindent $(1)$ $\pi-1$ lies in some prime $\mf{l}$ of $\bq(\W_{B_{\pi}})$ lying over $\ell$.

\noindent $(2)$ $\ell$ divides the order $\big|B_{\pi}(\bFq)\big|$.\end{Lem}

\begin{prf} $(1)\Rightarrow (2)$: If $\pi-1 \in \mf{l}$, then $\big|B_{\pi}(\bFq)\big| = \Nm(1-\pi)\in \mf{l}$ and since $\big|B_{\pi}(\bFq)\big|$ is a rational integer, $\big|B_{\pi}(\bFq)\big|\in \mf{l}\cap \bz = \ell\bz$.

\noindent $(2)\Rightarrow (1)$: Fix a prime $\mf{l}$ of $\bq(\W_{B_{\pi}})$ lying over $\ell$. If $\Nm(1-\pi) \in \ell\bz$, then there exists $\sigma \in \absq$ such that $1-\sigma^{-1}(\pi)\in \mf{l}$. Hence, $1-\pi\in \sigma(\mf{l})$.\end{prf}

\begin{Lem} The following are equivalent:

\noindent $(1)$ $\ell$ divides $\big|B_{\pi}(\bFq)\big|$ and $B_{\pi}$ has embedding degree $N$ with respect to $\ell$. \vspace{0.5mm}

\noindent $(2)$ There exists a prime $\mf{l}$ of $\bq(W_{B_{\pi}})$ lying over $\ell$ such that $\pi\equiv 1\pmod{\mf{l}}$ and $\Phi_N(\ov{\pi})\in \mf{l}$. \end{Lem}

\begin{prf} Let $\Phi_N(X)$ denote the $N$-th cyclotomic polynomial.

\noindent $(1)\Rightarrow (2)$: Since $\ell$ divides $P_{B_{\pi}}(1)$, it is clear that $\pi\equiv 1\pmod{\mf{l}}$ for some prime $\mf{l}$ lying over $\ell$. Now, $q\equiv \ov{\pi}\pmod{\mf{l}}$ and since $\Phi_N(q)\in \mf{l}$, it follows that $\Phi_N(\ov{\pi})\in \mf{l}$. \vspace{0.5mm}

\noindent $(2)\Rightarrow (1)$: If $\pi-1\in\mf{l}$, then $\Nm_{\bq(\pi)/\bq}(1-\pi)\in \Nm_{\bq(\pi)/\bq}(\mf{l}) = \ell\bz$ and hence, $\ell$ divides $\big|B_{\pi}(\bFq)\big|$. Furthermore, since $\Phi_N(\ov{\pi})\equiv \Phi_N(q) \pmod{\mf{l}}$, it follows that $\Phi_N(q)\in \mf{l}\cap \bz = \ell\bz$.\end{prf}

\begin{Lem} Let $B_{\pi}$ be an abelian variety over $\bFq$ and let $\ell$ be any rational prime not dividing $q$. Let $K/\bq$ be any number field and let $\ell\mc{O}_K = \mf{l}_1^{e_1}\cdots\mf{l}_r^{e_r}$ be the prime decomposition of $\ell$ in $K$. For each index $i$, let $d_i$ be the order of $q$ in $(\mc{O}_K/\mf{l}_i^{e_i})^*$. The integer $d_i$ is independent of the index $i$ and is the embedding degree of $B_{\pi}$ with respect to $\ell$. \end{Lem}

\begin{prf} For any index $i$ and integer $d$, $\Phi_d(q)\in \mf{l}_i$ if and only if $\Phi_d(q)\in \mf{l}_i\cap \bz = \ell\bz$. \end{prf}

The following well-known theorem shows that for elliptic curves, the embedding degree coincides with the the degree of the field extension generated by the $\ell$-torsion points.

\begin{Thm} $\mr{(Koblitz)}$ Let $E$ be an elliptic curve over $\bFq$ and let $\ell$ be a prime dividing $|E(\bFq)|$. Suppose $\ell\nmid (q-1)$. Then the embedding degree with respect to $\ell$ is $[\bFq(E[\ell]):\bFq]$.\end{Thm}

\begin{Lem} $\mr{([Kow03])}$ Let $B_{\pi}$ be a simple abelian variety over $\bFq$ and let $n$ be an integer prime to $q$. Then the degree $[\bFq(B_{\pi}[n]):\bFq]$ is given by \vspace{-0.15cm}$$[\bFq(B_{\pi}[n]):\bFq] = \mr{inf}\{d\geq 1: \pi^d\equiv 1\Mod{n} \text{ in } \End(A)\}. $$\end{Lem}

\noindent The following corollary follows immediately.

\begin{Corr} Let $\ell$ be a prime that does not divide $q$ and let $\mf{l}$ be a prime of $\wti{L}:=\bq(\W_{B_{\pi}})$ lying over $\ell$. For every $\sigma\in\absq$, let $d_{\sigma}$ be the order of $\sigma(\pi)$ in $(\mc{O}_{L}/\mf{l}^e\mc{O}_{L})^*$ where $e$ is the ramification index of $\ell$ in $\bq(\W_{B_{\pi}})/\bq$ . Then the degree $[\bFq(B_{\pi}[\ell]):\bFq]$ is given by $\mr{lcm}(d_{\sigma})_{\sigma}$.\end{Corr}

\begin{prf} We write $d:=\mr{lcm}(d_{\sigma})_{\sigma}$. Let $\mf{l}_1$ be any prime of $\bq(\W_{B_{\pi}})$ lying over $\ell$. Then there exists $\sigma\in \absq$ such that $\mf{l} = \sigma(\mf{l}_1)$. So $\pi^{d_{\sigma}} \equiv 1 \pmod{\mf{l}_1^e}$ and hence, $\pi^d\equiv 1\pmod{\ell}$.

Conversely, if $\pi^{d_1}\equiv 1\pmod{\ell}$ for some integer $d_1$, then in particular, $\pi^{d_1}\equiv 1\pmod{\sigma^{-1}(\mf{l}^e)}$ and hence, $\sigma(\pi)^{d_1}\equiv 1\pmod{\mf{l}^e}$. Hence, $d$ divides $d_1$.\end{prf}

\noindent \textbf{Remark.} The assumption that $\gcd(n,q) = 1$ is necessary to ensure that $B_{\pi}[n]$ is etale. Clearly, if $\pi^d \equiv 1\pmod {\ell\mc{O}_{\bq(\pi)}}$, then $\ov{\pi}^d \equiv 1\pmod {\ell\mc{O}_{\bq(\pi)}}$ and hence, $q^d\equiv 1\pmod{\ell}$ since complex conjugation preserves the ideal $\ell\mc{O}_{\bq(\pi)}$. So the embedding degree $n$ always divides the degree $[\bFq(B_{\pi}[\ell]):\bFq]$ of the field generated by all $\ell$-torsion points. Thus, we have the inclusion \vs $$\bFq\sub \bF _{q^n}\sub \bFq(B_{\pi}[\ell]) \sub \bF _{q^{\ell-1}}.  $$ But unlike the case of elliptic curves, the degree $[\bFq(B_{\pi}[\ell]):\bFq]$ can be substantially larger than the embedding degree in the case of higher dimensional abelian varieties.

In pairing based cryptography, it is eminently desirable for the prime $\ell$ to be substantially larger than the embedding degree with respect to $\ell$. One way to achieve this is to construct an abelian variety $B$ over a finite field $\bFq$ such that the \textit{full embedding degree}, i.e. the degree of the extension generated by the $\ell$-torsion points is substantially larger than the embedding degree. By Koblitz's theorem, this is not possible for elliptic curves since the two notions coincide. However, as shown in the next section, this is possible for dimension two or greater.

\subsection{\fontsize{11}{11} Newton Polygons}


Let $B$ be an abelian variety over an algebraically closed field $k$ of characteristic $p>0$. The group scheme $B[p^{\infty}]$ is a $p$-divisible group of rank $\leq \dim B$. Let $D(B[p^{\infty}])$ denote the Dieudonne module. Then $D(B[p^{\infty}])\otimes_k W(k)[\frac{1}{p}]$ is a direct sum of pure isocrystals by the Dieudonne-Manin classification theorem. Let $\lamb_1<\cdots<\lamb_r$ be the distinct slopes and let $m_i$ denote the multiplicity of $\lamb_i$. The sequence $m_1\times \lamb_1,\cdots, m_r\times \lamb_r$ is called the \textit{Newton polygon} of $B$.

\begin{Def} \normalfont A Newton polygon is \textit{admissible} if it fulfills the following conditions:

\noindent 1. The breakpoints are integral, meaning that for any slope $\lambda$ of multiplicity $m_{\lambda}$, we have $m_{\lambda}\cdot \lambda \in \bz$.
 
\noindent 2. The polygon is \textit{symmetric}, meaning that for each slope $\lambda$, the slopes $\lambda$ and $1-\lambda$ have the same multiplicity.\end{Def}

Let $\pi$ be a Weil $q$-integer and let $B_{\pi}$ be the corresponding simple abelian variety over $\bFq$. Then the Newton slopes of $B_{\pi}$ are given by $\{v(\pi)/v(q)\}_v$ where $v$ runs through the places of $\bq(\pi)$ lying over $p$. In particular, the Newton polygon is symmetric and hence, all slopes lie in the interval $[0,1]$. A far more subtle fact is that the converse is also true, i.e. every admissible Newton polygon occurs as the Newton polygon of a simple abelian variety over a finite field. This was formerly known as Manin's conjecture until proven by Serre. We refer the reader to [Tat69] for the proof.

\begin{Thm} The Newton polygon of an abelian variety over a finite field is admissible. Conversely, any admissible polygon occurs as the Newton polygon of some abelian variety over a finite field of any prescribed characteristic.\end{Thm}

\begin{Def} \normalfont An simple abelian variety of dimension $g$ over a field $F$ is \textit{ordinary} if its Newton polygon is $g\times 0\;,\;g\times 1$. 
\end{Def}

For an abelian variety corresponding to a Weil $q$-integer $\pi$, the Newton slopes are given by $\frac{v(\pi)}{v(q)}$ as $v$ runs through the places of $\bq(\pi)$ lying over the prime $p$ that divides $q$. The element $\frac{v(\pi)}{v(q)}\in \bq/\bz$ is $0$ precisely when $v$ does not divide $\pi$. Similarly, $\frac{v(\pi)}{v(q)} = 1$ if and only if $v$ divides $\pi$ but not $\ov{\pi}$. Thus, $B_{\pi}$ being ordinary is equivalent to $\pi$ being coprime to $\ov{\pi}$ in $\bq(\pi)$ (or in any field extension of $\bq(\pi)$).

\begin{Def} \normalfont An simple abelian variety of dimension $g$ over a field $F$ is \textit{supersingular} if its Newton polygon is $2g\times 1/2$. \end{Def}

\section{\fontsize{11}{11}\selectfont The supersingular case}

We first show that the embedding degree can be arbitrarily large. The supersingular case will be key here.

\begin{Def} \normalfont An abelian variety $A$ over a finite field of characteristic $p$ is supersingular if all of its Newton slopes are $\frac{1}{2}$. \end{Def}

For dimension $1$ or $2$, one may alternatively define a supersingular abelian variety as one with no $p$-torsion points or equivalently with no $0$ as a Newton slope. But that does not extend to higher dimensions, as evidenced by absolutely simple abelian threefolds with Newton polygon $3\times \frac{1}{3}\;,\; 3\times \frac{2}{3}$. It is a well-known fact that a $g$-dimensional supersingular abelian variety is isogenous over the algebraic closure to the $g$-th power of the supersingular elliptic curve. 

\begin{Lem} Let $\pi$ be a Weil $q$-integer corresponding to a simple abelian variety of dimension $g$. The following are equivalent:\\
\noindent $1$. $B_{\pi}$ is supersingular\\
\noindent $2$. $\pi^N\in \bz$ for some integer $N$.\\
\noindent $3$. $B_{\pi}\times _{\bFq} \ov{\bF}_q$ is isogenous to a power of the supersingular elliptic curve.\\
\noindent $4$. $\End^0(B_{\pi})\cong \Mat_{g}(\bq_{p,\infty})$, the quaternion algebra over $\bq$ ramified only at $p$ and $\infty$\end{Lem}

\begin{Lem} Let $\pi_{_1}, \pi_{_2}$ be Weil $q$-integers. The corresponding simple abelian varieties over $\bFq$ are twists of each other if and only if $\sigma(\pi_{_1})\cdot \pi_{_2}^{-1}$ is a root of unity for some $\sigma\in \absq$.\end{Lem}
\begin{prf} ($\Rightarrow$): Suppose $B_{\pi_{_1}}$, $B_{\pi_{_2}}$ are isogenous over some extension $\bF _{q^N}$. Then $\pi_{_1}^N$ and $\pi_{_2}^N$ are conjugates over $\bq$. Choose $\sigma \in \absq$ such that $\sigma(\pi_{_1}^N) = \pi_{_2}^N$. Then $(\sigma(\pi_{_1})\cdot \pi_{_2}^{-1})^N = 1$ and hence, $\sigma(\pi_{_1})\cdot\pi_{_2}^{-1}$ is an $N$-th root of unity.\vspace{1mm}

\noindent $(\Leftarrow)$: Suppose $(\sigma(\pi_{_1})\cdot \pi_{_2}^{-1})^N = 1$ for some $\sigma\in \absq$. Then $\sigma(\pi_{_1}^N) = \pi_{_2}^N$ and hence, $B_{\pi_{_1}^N}$ is isogenous to $B_{\pi_{_2}^N}$.\end{prf}

In particular, if $\pi$ is a Weil $q$-integer corresponding to a supersingular abelian variety, then $\frac{\pi}{\sqrt{q}}$ is a root of unity.

\vspace{1.5mm}

\begin{Corr} Let $B$ be a simple supersingular abelian variety over a finite field $\bFq$. Then $\dim B = \phi(N)$ or $\frac{\phi(N)}{2}$ for some integer $N$. \end{Corr}

\begin{prf} Let $\pi$ be the corresponding Weil $q$-integer. Then $\pi = \sqrt{q}\cdot \zeta_{_N}$ where $N$ is an integer and $\zeta_{_N}$ is a primitive $N$-th root of unity. For brevity, write $D_{\pi}:=\End^0(B)$ and let $p$ be the characteristic of $\bFq$. By Honda-Tate theory, $2\dim B = [\bq(\pi):\bq]\cdot [D_{\pi}:\bq(\pi)]^{1/2}$. Now, $$[\bq(\pi):\bq] = \begin{cases}
\phi(N) & \text{if } \sqrt{q}\in \bq(\zeta_{_N})\\
2\phi(N) & \text{if } \sqrt{q}\notin \bq(\zeta_{_N}).\end{cases}$$ 

Furthermore, the Newton slopes are all $\frac{1}{2}$. So the Hasse invariants of $D_{\pi}$ are $\frac{1}{2}\cdot [\bq(\pi)_v:\bq_p]$ at the places  $v$ lying over $p$. Thus, the least common denominator of the Hasse invariants is either $1$ or $2$ and hence, $[D_{\pi}:\bq(\pi)] = 1$ or $4$. 

If $\sqrt{q}\notin \bq(\zeta_{_N})$, then $p$ is ramified in $\bq(\pi)/\bq$ and hence, the Hasse invariants of $D_{\pi}$ are $0$ in $\bq/\bz$, meaning that $D_{\pi} = \bq(\pi)$. This completes the proof.\end{prf}

\vspace{1mm}

\begin{Prop} Let $N$ be any integer and let $q$ be a prime power that is a perfect square. Then there exists a simple supersingular abelian variety $B$ over $\bFq$ and a rational prime $\ell$ such that:

\noindent -the group $B(\bFq)$ of $\bFq$-points has order divisible by $\ell$.\\ 
\noindent -the embedding degree of $B$ with respect to $\ell$ is $N$.\end{Prop}

\begin{prf} Let $\pi$ be the Weil $q$-integer $\sqrt{q}\cdot \zeta_{_{2N}}$ where $\zeta_{_{2N}}$ is a primitive $2N$-th root of unity. Let $B$ be the corresponding simple abelian variety over $\bFq$. Then $B_{\pi}$ is supersingular. 

Now, the order of $B(\bFq)$ is given by: 
\begin{equation*}
\label{eq:pareto mle2}
 \begin{aligned} 
P_{B}(1) = \prod\limits_{\sigma\in \GalWB}(1-\sigma(\pi))\\
 = \prod\limits_{\gcd(i,\;2N)=1}(1-\sqrt{q}\zeta_{_{2N}}^i)\\ = \prod\limits_{\gcd(i,\;2N)=1}(\zeta_{_{2N}}^{-i}-\sqrt{q})\\ = \Phi_{2N}(\sqrt{q}).
\end{aligned}
\end{equation*}
Thus, any prime $\ell$ dividing $\Phi_{2N}(\sqrt{q})$ but not dividing $2N$ fulfills the hypothesis of the proposition. Since such a prime always exists unless $(q, N) = (4, 3)$, this completes the proof.\end{prf}

Clearly, the prime $\ell$ constructed in the last proposition is $\equiv 1\pmod N$. So choosing a moderately large value for $N$ results in $\ell$ being a substantially large prime.

The simple abelian variety constructed in the last proposition is of dimension ${\phi(2N)}{}$ and is far from absolutely simple. Later, we shall provide a construction that yields an absolutely simple abelian variety of any prescribed dimension $g$ with embedding degree any prescribed integer $N$ with respect to some prime $\ell$. 

\begin{Prop} Let $B$ be a simple supersingular abelian variety of dimension $g$ over $\bFq$. Then for any prime $\ell$ dividing $|B(\bFq)|$, the embedding degree of $B$ with respect to $\ell$ is $\leq 4g^2$ \end{Prop}

\begin{prf} Let $\pi$ be the associated Weil number and let $p$ be the prime dividing $q$. Then $\pi = \sqrt{q}\cdot\zeta_{_N}$ where $\zeta_{_N}$ is a primitive $N$-th root of unity for some integer $N$. The case where $N \in \{1,2\}$ is trivial. So we may assume $N\geq 3$. Thus, $\bq(\pi) = \bq(\zeta_N)$, which is a CM field.

Now, $$2g = [\bq(\pi):\bq]\cdot [D_{\pi}:\bq(\pi)]^{1/2} = \phi(N)\cdot [D_{\pi}:\bq(\pi)] . $$ Furthermore, $v(\pi)/v(q) = 1/2$ at the places of $\bq(\zeta_N)$ lying over $p$ and $0$ elsewhere. Thus, $D_{\pi}$ is either $\bq(\zeta_N)$ or the quaternion algebra over $\bq(\zeta_N)$ ramified exclusively at the primes lying over $p$, depending on whether the local degree of a prime of $\bq(\zeta_N)$ lying over $p$ is even or odd. Hence, $g  = \phi(N)$ or $\frac{1}{2}\cdot \phi(N)$. 

Now, $\pi\equiv 1\pmod{\mf{l}}$ for some prime $\mf{l}$ of $\bq(\W_B)$ lying over $\ell$. So $q\equiv \zeta_{_N}^{-2}\pmod{\mf{l}}$ and hence, $\ell$ divides $\Phi_M(q)$ where $M = \frac{N}{\gcd(2,N)}$. Thus, the embedding degree divides $N$. Now, $\phi(N) = g$ or $2g$ and hence, $N \leq 4g^2$.\end{prf}

\noindent \textbf{Remark:} This is just a straightforward generalization of the observation that supersingular elliptic curves have embedding degrees at most $6$. This arises from the fact that imaginary quadratic fields do not contain any roots of unity other than those of torsion $4$ or $6$. The upshot is that supersingular abelian varieties make largely unsuitable candidates for pairing based cryptography.

\section{\fontsize{11}{11}\selectfont Base changes with prescribed embedding and full embedding degrees}

The following proposition demonstrates that Koblitz's aforementioned theorem does not generalize to higher dimensions. In fact, if an abelian variety of dimension $\geq 2$ satisfies some rather generic conditions, there exists a base change with prescribed embedding and full embedding degrees.

\begin{Thm}\label{generalize} Let $B_{\pi}$ be an absolutely simple abelian variety of dimension $\geq 2$ over a finite field $\bFq$ such that the multiplicative group $\Phi_{B_{\pi}}$ generated by $\W_{B_{\pi}}$ is free abelian of rank $\geq 3$. Let $m$, $n$ be prescribed integers. Then there exists a base change $B_{\wti{\pi}}$ of $B_{\pi}$ to some finite extension $\bF _{\wti{q}}$ and a prime $\ell$ such that:

\noindent - $\ell$ divides $|B_{\wti{\pi}}(\bFq)|$\\
\noindent - $B_{\wti{\pi}}$ has embedding degree $n$ with respect to $\ell$.\\
\noindent - $B_{\wti{\pi}}$ has full embedding degree $m\cdot n$ with respect to $\ell$, meaning $[\bF _{\wti{q}}(B_{\wti{\pi}}[\ell]):\bF _{\wti{q}}] = m\cdot n$. \end{Thm}

\begin{prf} We write $L := \bq(\pi)$ and denote its Galois closure over $\bq$ by $\wti{L}$. Choose a conjugate $\pi_{_1}$ of $\pi$ such that $\pi_{_1}$ lies outside the subgroup of $\Phi_{B_{\pi}}$ generated by $\pi$ and $\ov{\pi}$. By Kummer theory, the field $\wti{L}(\zeta_{mn},\sqrt[mn]{\pi},\sqrt[mn]{\ov{\pi}},\sqrt[mn]{\pi_{_1}})$ is abelian over $\wti{L}(\zeta_{mn})$ with Galois group $(\bz/mn\bz)^3$. We choose a prime $\mf{l}$ of $\wti{L}(\zeta_{mn})$ such that:

\noindent - $\mf{l}$ is of local degree one over $\bq$.

\noindent - $\mf{l}$ splits completely in $\wti{L}(\zeta_{mn},\sqrt[m\cdot n]{\pi})$.

\noindent - The  decomposition field of $\mf{l}$ in $\wti{L}(\zeta_{mn},\sqrt[m\cdot n]{\ov{\pi}})$ is $\wti{L}(\zeta_{mn},\sqrt[m]{\ov{\pi}}).$

\noindent - $\mf{l}$ is inert in $\wti{L}(\zeta_{mn},\sqrt[{m\cdot n}]{\pi_{_1}})$.

Note that by the Chebotarev density theorem, the set of primes fulfilling these conditions has positive density. Let $\ell$ be the rational prime lying under $\mf{l}$. By construction, $\ell\equiv 1\pmod{m\cdot n}$ and \vs $$ \pi^{\frac{\ell-1}{m\cdot n}}\equiv 1\Mod{\mfl} \;,\;  \ov{\pi}^{\frac{\ell-1}{m\cdot n}}\equiv \zetamn^{m\cdot i}\Mod{\mfl}\;,\; \pi_{_1}^{\frac{\ell-1}{m\cdot n}}\equiv \zetamn^j\Mod{\mfl} $$ for some integers $i$, $j$ co-prime to $m\cdot n$. We write \vspace{-1mm}$$\wti{\pi}:=\pi^{^{\frac{\ell-1}{m\cdot n}}}\;\;,\;\;\wti{q}:=q^{^{\frac{\ell-1}{m\cdot n}}}\vspace{-1mm}$$ for brevity. Now, $\wti{\pi}$ is a Weil $\wti{q}$-integer corresponding to the base change of $B_{\pi}$ to $\bF _{\wti{q}}$. Furthermore, $\wti{\pi}\equiv 1\pmod{\mf{l}}$ and hence, $\ell$ divides $|B_{\pi}(\bF _{\wti{q}})|$. On the other hand, the complex conjugate $\ov{\wti{\pi}}$ of $\wti{\pi}$ has order $n$ in $(\mc{O}_{\wti{L}(\zeta_{_N})}/\mf{l})^*$ and since $\wti{q}\equiv \ov{\wti{\pi}}\pmod{\mf{l}}$, it follows that $\Phi_n(\wti{q})\in \mf{l}$. Thus, $B_{\wti{\pi}}$ has embedding degree $n$ with respect to $\ell$.

Now, the full embedding degree $[\bF _{\wti{q}}(B_{\wti{\pi}}[\ell]):\bF _{\wti{q}}]$ of $B_{\wti{\pi}}$ with respect to $\ell$ is is given by \vspace{-1mm}$$\mr{inf}\{d\geq 1: \sigma(\wti{\pi})^d\equiv 1\Mod{\mf{l}} \;\;\forall\;\;\sigma\in\absq \}. \vspace{-1mm} $$ Since $\ell$ splits completely in $\wti{L}(\zetamn)$, the abelian group $(\mc{O}_{\wti{L}(\zetamn)}/\ell\mc{O}_{\wti{L}(\zetamn)})^*$ is a product of $[{\wti{L}(\zetamn)}:\bq]$ copies of $\bFl^*$. In particular, the group is $(\ell-1)$-torsion and hence, the order of any $\absq$-conjugate of $\wti{\pi}$ in $(\mc{O}_{{\wti{L}(\zetamn)}}/\ell\mc{O}_{{\wti{L}(\zetamn)}})^*$ divides $m\cdot n$. 

By construction, $\wti{\pi}_{_1}:=\pi_{_1}^{^{\frac{\ell-1}{m\cdot n}}}$ has order $m\cdot n$ in $(\mc{O}_{\wti{L}(\zeta_{_{mn}})}/\mf{l})^*$. Thus, $m\cdot n$ divides the degree $[\bF _{\wti{q}}(B_{\wti{\pi}}[\ell]):\bF _{\wti{q}}]$, which forces equality. \end{prf}

\noindent \textbf{Remark.} Note that the condition that $\Phi_{B_{\pi}}$ is free abelian of rank $\geq 3$  is quite generic when the dimension is $\geq 2$. The only class of exceptions is that of the abelian varieties of type $\IV(1,d)$ which we will encounter in the next section.

\begin{Prop} \label{rank} Let $B_{\pi}$ be an absolutely simple abelian variety corresponding to a Weil $q$-integer $\pi$ such that $B_{\pi}$ has dimension $\geq 2$ and $\Phi_{B_{\pi}}$ is free abelian. Unless $\big|W_ {B_{\pi}}\big|\leq 2$,  the abelian group $\Phi_{B_{\pi}}$ has rank $\geq 3$. \end{Prop}

\begin{prf} Since any supersingular abelian variety decomposes as a product of supersingular elliptic curves over the algebraic closure, we may assume $B_{\pi}$ is not supersingular. Thus, there do not exist non-zero integers $x,y$ such that $\pi^x = q^y$.

Suppose $\Phi_{B_{\pi}}$ has rank $\leq 2$. Let $\sigma \in \absq$ be an automorphism such that $\sigma(\pi)\notin \{\pi,\ov{\pi} \}$. Then $\sigma(\pi) = \pi^{k_1}\cdot q^k$ for some non-zero integers $k,k_1$ such that $k_1+2k = 1$, $(k,k_1)\neq (1,-1)$. So \vs $$ \sigma^{j}(\pi) = \pi^{k_1^j}\cdot q^{k(k_1^{j}-1)/(k_1-1)} .  $$  In particular, if $\sigma$ has order $n$, we have $\pi^{k_1^n} = q^{-k(k_1^{n}-1)/(k_1-1)}$, a contradiction since $B_{\pi}$ is not supersingular. \end{prf}

\vspace{1mm}

\noindent We note that the condition that $\Phi_{B_{\pi}}$ is free abelian is necessary for the preceding proposition to work. We describe a class of $3$-dimensional absolutely simple abelian varieties $B_{\pi}$ for which $\big|\W_ {B_{\pi}}\big| = 6$ and $\Phi_{B_{\pi}}$ is of rank $2$.

Let $p\neq 3$ be a prime and let $\pi_1$ be an ordinary Weil $p$-integer corresponding to an ordinary elliptic curve over $\bFp$. Then $\bq(\pi_1)$ is an imaginary quadratic field in which $p$ splits into two distinct principal prime ideals $(\pi_1)$ and $(\ov{\pi}_1)$. Let $\pi$ be a zero of $X^3- \pi_1^{2}\cdot\ov{\pi}_1$. Then $\pi$ is a Weil $p$-integer corresponding to a simple abelian variety $B_{\pi}$ (up to isogeny) of dimension $3$ over $\bFp$. The set $\W_ {B_{\pi}}$ is given by \vs $$ \W_ {B_{\pi}} = \{\pi\cdot \zeta_3^{i}\;,\; \ov{\pi}\cdot \zeta_3^{i}:\;\;i=0,1,2  \} $$ where $\zeta_3$ is a cube root of unity. Thus, $\big|\W_ {B_{\pi}}\big| = 6$ and $\Phi_{B_{\pi}}$ has rank $2$. $B_{\pi}$ has Newton polygon \vs $$ 3 \times \frac{1}{3}\;\;,\;\; 3 \times \frac{2}{3}$$ and its endomorphism algebra is isomorphic to the degree six CM field $\bq(\pi)$.

A base change of $B_{\pi}$ to the field $\bF _{p^3}$ yields the abelian variety $B_{\pi^3}:= B_{\pi}\times_{\bFp} \bF _{p^3} $ with \vs $$\W _ {B_{\pi^3}} = \{\pi_1^2\cdot\ov{\pi}_1\;,\;\pi_1\cdot\ov{\pi}_1^2\}\;\;\,,\;\;\Phi_{B_{\pi^3}} = \langle \pi_1,\ov{\pi}_1 \rangle $$ and in particular, $\Phi_{B_{\pi^3}}$ is free abelian of rank $2$. The abelian variety $B_{\pi^3}$ is absolutely simple and its endomorphism algebra $D_{\pi^3}$ is a $9$-dimensional central division algebra over the imaginary quadratic field $\bq(\pi_1)$ with Hasse invariants $\frac{1}{3}$ and $\frac{2}{3}$ at the two places lying over $p$. The degree six CM field $\bq(\pi)$ has an embedding in this division algebra.

We thank the anonymous referee for recommending the inclusion of this counter-example in the case where $\Phi_{B_{\pi}}$ is not assumed to be free abelian.

\section{\fontsize{11}{11}\selectfont The ordinary case}

Having seen why supersingular abelian varieties are unsuitable candidates, we now look at the other end of the spectrum, which is ordinary abelian varieties. The primary task of this subsection is to address the following: \vspace{2mm}

\noindent \textbf{Problem:} For prescribed integers $g$, $N$, construct an ordinary $g$-dimensional abelian variety $B$ over a finite field $\bFq$ with embedding degree $N$ with respect to a prime $\ell$ dividing $|B(\bFq)|$. \vspace{2mm}

Clearly, for this to be possible, the prime $\ell$ must split completely in $\bq(\zeta_{_N})$. We first observe that such a construction is not always possible if the prime $\ell\equiv 1\pmod{N}$ is prescribed. \vspace{1mm}

\noindent \textbf{Example.} Let $L$ be a CM field of degree $2g$ such that $L/\bq$ is cyclic and $L$ is linearly disjoint with $\bq(\zeta_{_N})$. Let $\ell$ be a rational prime that splits completely in $\bq(\zeta_{_N})$ but is inert in $L/\bq$. 

Let $\pi$ be an ordinary Weil number in $L$ corresponding to a simple abelian variety over the finite field $\bFq$ of size $q := \pi\cdot\ov{\pi}$. If $\ell$ divides $\big|B_{\pi}(\bFq)\big|$, then $\pi\equiv 1\pmod{\mf{l}}$ for some prime $\mf{l}$ of $\ell$ lying over $\ell$. But since $\ell\mc{O}_L$ is assumed to be a prime ideal, it follows that $\pi\equiv 1\pmod{\ell}$. So $\ov{\pi}\equiv 1 \pmod{\ell}$ and hence, $q\equiv 1\pmod{\ell}$. Thus, any Weil number in $L$ corresponds to an abelian variety with embedding degree $1$ with respect to $\ell$. \vspace{2mm}

\noindent \textbf{The motivation.} Let $B_{\pi}$ be an ordinary abelian variety of dimension $g$ over a finite field $\bFq$ of characteristic $p$. Write $q = p^k$. Suppose the number field $L:=\bq(\pi)$ is abelian. Then $p$ splits completely in $L$. Now, the principal ideals $\pi\mc{O}_L$ and $\ov{\pi}\mc{O}_L$ have prime decompositions \vs $$\pi\mc{O}_L = \prod\limits_{i=1}^g \mfp_i\;\;\; ,\;\;\;\ov{\pi}\mc{O}_L = \prod\limits_{i=1}^g \ov{\mfp}_i$$ where the $\mfp_1,\; \ov{\mfp}_1,\;\cdots,\; \mfp_g,\; \ov{\mfp}_g$ are distinct primes. Let $\sigma_i\in \Gal(L/\bq)$ be the element such that $\sigma_i(\mfp_1) = \mfp_i$. Then $\Phi:=\{\sigma_i:1\leq i\leq g \}$ is a CM type for $L$.

Thus, when it comes to constructing ordinary abelian varieties over prime fields, CM types are a natural place to look.

\vspace{2mm}

\noindent \textbf{The reflex field.} Let $L$ be a CM field of degree $2g$ and let $\wti{L}$ be its Galois closure over $\bq$. A \textit{CM-type} $\Phi = \{\phi_1,\cdots,\phi_g\}$ for $L$ is a set of embeddings $L\hra \wti{L}$ such that $\ov{\Phi} = \{\ov{\phi}_1,\cdots,\ov{\phi}_g\}$ is disjoint from $\Phi$ and the union $\Phi\coprod \ov{\Phi}$ is the complete set of embeddings $L\hra \wti{L}$. The map $$\mr{Nm}_{\Phi}:L\lra \wti{L}\;\;,\;\; x\mapsto \prod\limits_{\phi\in\Phi}\phi(x) $$ is called the \textit{type norm} with repect to $\Phi$. Clearly, this map satisfies $$\mr{Nm}_{\Phi}(x)\cdot {\mr{Nm}_{\ov{\Phi}}(x)} = \mr{Nm}_{L/\bq}(x).$$   

Let $\wti{L}$ denote the Galois closure of $L$ over $\bq$ and let $S$ be the union of left cosets of $\Gal({L}/\bq)$ in $\Gal(\wti{L}/\bq)$. Set $\Psi:= \{\sigma^{-1}:\sigma \in S\}$ and let $\what{H}$ be the stabilizer of $\Psi$ in $\Gal(\wti{L}/\bq)$. Then $\what{H}$ defines a subfield of $\wti{L}$ and we have $\what{H} = \{\gamma\in \Gal(\wti{L}/\bq): \Psi\gamma = \Psi \}$. These cosets define an embedding of $\what{L}$ into $\wti{L}$. We call $(\what{L})$ the reflex field of $(L,\Phi)$ and $\Psi$ the reflex type.

The map \vspace{-0.15cm}$$\Nm_{\Psi}:\wti{L}\lra \wti{L}\;\;\; ,\;\;\; x\mapsto \prod\limits_{\psi\in\Psi}\psi(x)$$ maps the field $\what{L}$ to $L$. This is called the \textbf{reflex norm} of $\Psi$. In the generic case, the field $\wti{L}$ is of degree $2^g \cdot g!$ and $\what{L}$ is of degree $2^g$ with $\Gal(\wti{L}/\what{L}) = S_g$, the symmetric group on $g$ letters. 

\bigskip

\begin{Thm} Let $L$ be a fixed CM field of degree $\geq 4$ and let $m$, $n$ be prescribed integers. Let $p$ be a rational prime that splits completely in the Hilbert class field $\mr{H}(L)$ of $L$. There exists an absolutely simple ordinary abelian variety $B$ over the prime field $\bFp$, a $p$-power $q$ and a prime $\ell$ such that:

\noindent - $\End^0(B) \cong L$ 

\noindent - $\ell$ divides $|B(\bFq)|$ and $B\times _{\bFp}\bFq$ has embedding degree $n$ with respect to $\ell$.

\noindent - $B\times _{\bFp}\bFq$ has full embedding degree $m\cdot n$ with respect to $\ell$.\end{Thm}

\begin{prf} Let $K$ be the maximum totally real subfield of $L$ and let $\wti{K}$, $\wti{L}$ be the respective Galois closures over $\bq$. We fix a primitive CM type $\Phi$ for $L$ and denote the reflex of $(L,\Phi)$ by $(\what{L}, \Psi)$. 

Note that the Hilbert class field $\mr{H}(\wti{L})$ of $\wti{L}$ coincides with the Galois closure $\wti{\mr{H}(L)}$ of $\mr{H}(L)$ over $\bq$. Since $p$ splits completely in $\mr{H}(L)$, it also splits completely in the Galois closure $\wti{\mr{H}(L)} = \mr{H}(\wti{L})$ of $\mr{H}(L)$. We choose a prime $\wti{\mf{p}}$ of $\wti{L}$ lying over $p$. Thus,

\noindent -$\wti{\mf{p}}$ has local degree one over $\bq$.\\
-$\wti{\mf{p}}$ splits completely in the Hilbert class field $\mr{H}(\wti{L})$.

The second condition ensures that $\wti{\mf{p}}$ is principal. We write $\wti{\mf{p}} = \wti{\al}\mc{O}_{\wti{L}}$. Define $\what{\mfp}:= \wti{\mf{p}}\cap \what{L}$ and $\what{\al}:= \Nm_{\wti{L}/\what{L}}(\wti{\al})$. Then $\what{\mfp}$ has local degree one over $\bq$ and $\what{\mfp} =\what{\al}\mc{O}_{\what{L}}$. In particular, the prime $\what{\mfp}$ of $\what{L}$ is principal. Set $\pi:= \Nm_{\Psi}(\what{\al})$. Then $\pi\cdot \ov{\pi} = p$ and since the map $\psi:\wti{L}\lra \wti{L}$ maps $\what{L}$ to $L$, it follows that $\pi\in \mc{O}_L$. Furthermore, the ideals $(\pi)$ and $(\ov{\pi})$ of $\mc{O}_{{L}}$ are coprime and hence, $\pi$ is an ordinary Weil $p$-integer. 

Now, the multiplicative group $\Phi_{B_{\pi}}$ generated by the Galois conjugates of $\pi$ is free abelian of rank $\geq 3$. Suppose, by way of contradiction, that $\Phi_{B_{\pi}}$ is of rank $\leq 2$. Let $\pi_1$ be a conjugate of $\pi$ other than $\pi$, $\ov{\pi}$. Then there exist non-zero integers $r_1$, $r_2$, $r_3$ such that $\pi_1^{r_1} = \pi^{r_2}\ov{\pi}^{r_3}$. If $r_2, r_3 >  0$, it follows that $\pi$ is divisible all primes of $\bq(\pi)$ lying over $p$, a contradiction since $B_{\pi}$ is ordinary. If $r_2<0$, it follows that $(\pi_1)$ divides $\ov{\pi}^{r_3}$ in $\mc{O}_{\bq(\pi)}$, a contradiction since $(\pi_1,\ov{\pi}_1) = (1)$ in $\bq(\pi_1)$.

The theorem now follows from Theorem \ref{generalize}. \end{prf}

\vspace{2mm}

\begin{Thm} Let $L$ be a fixed CM field of degree $2g\geq 4$. Let $m$, $n$ be prescribed integers and $\ell$ a rational prime $\equiv 1 \pmod{m\cdot n}$ that splits completely in $L/\bq$. There exists a simple ordinary abelian variety $B_{\pi}$ over a prime field $\bFp$ such that:

\noindent - $\End^0(B_{\pi})\cong L$.\\
- $\ell$ divides the order $\big|B_{\pi}(\bFp)\big|$ of the group of $\bFp$-points.\\
- $B_{\pi}$ has embedding degree $n$ with respect to $\ell$.\\
- $B_{\pi}$ has full embedding degree $m\cdot n$ with respect to $\ell$, meaning $[\bFp(B_{\pi}[\ell]):\bFp] = m\cdot n$.\end{Thm}

\begin{prf} Let $\wti{L}$ be a Galois closure of $L$ over $\bq$. Then $\wti{L}$ is a Galois CM field. Since $\ell$ splits completely in $L$ and is $\equiv 1\Mod{m\cdot n}$, it splits completely in $\wti{L}(\zetamn)$.

We fix a primitive CM type $\Phi$ for $\wti{L}$ that is stable under the action of $\Gal(\wti{L}/L)$. So $\Gal(\wti{L}/\bq) = \Phi \cup \ov{\Phi}$ and the image of the map \vspace{-1mm} $$ \mr{Nm}_{\Phi}: \wti{L}\lra \wti{L}    \vspace{-1mm} $$ lands in $L$. We may assume without loss of generality that $\Phi$ is not preserved or swapped for $\ov{\Phi}$ under the action of $\absq$, i.e. there exists a $\tau\in \absq$ and $\sigma_1,\sigma_2\in \Phi$ such that $\sigma_1,\ov{\sigma}_2\in \tau(\Phi)$. 

Let $\mfl$ be a prime of $\wti{L}$ lying over $\ell$. Then \vs $$\ell\mc{O}_{\wti{L}} =\prod\limits_{\phi\in \Phi} \phi(\mfl)\ov{\phi}(\mfl) .\vs $$ Let $r\in \bz$ be a rational integer such that $\Phi_{mn}(r)\in \mfl$. Then $\Phi_{mn}(r)\in \ell\bz$. 

Let $\wti{L}^{\ell \mc{O}_{\wti{L}}}$ denote the ray class field of the modulus $\ell\mc{O}_{\wti{L}}$ and let $\wti{L}^{\sigma(\mfl)}$ denote the ray class field of the modulus $\sigma(\mfl)$ for any $\sigma\in \absq$. Then $$\Gal(\wti{L}^{\ell \mc{O}_{\wti{L}}}/\wti{L})\cong \prod\limits_{\sigma \in \Gal(\wti{L}/\bq)} \Gal(\wti{L}^{\sigma(\mfl)}/\wti{L}) $$ and each $\Gal(\wti{L}^{\sigma(\mfl)}/\wti{L})$ is isomorphic to the quotient group $\mc{O}_L^{\times}/(\mfl + 1)$, which is cylic of order $\ell-1$ has a subgroup of order $mn$. The Chebotarev density theorem applied to the extension $\wti{L}^{\ell \mc{O}_{\wti{L}}}/\wti{L}$ implies that the Dirichlet density of principal prime ideals of $\wti{L}$ generated by elements $\al$ satisfying the following properties is positive:

\noindent - $\al$ generates a principal prime ideal in $\mc{O}_{\wti{L}}$.

\noindent - $\al\equiv r\Mod{\sigma_1^{-1}\mfl}\;$, $\;\al\equiv r^{m-1}\Mod{\ov{\sigma}_2^{-1}(\mfl)}$.

\noindent - $\al\equiv 1\Mod{\sigma(\mfl)}$ for every $\sigma\in \Gal(\wti{L}/\bq)\setminus \{\sigma_1^{-1},\ov{\sigma}_2^{-1} \}$

Let $p$ be the rational prime lying under $\al\mc{O}_{\wti{L}}$. We set $\pi:= \mr{Nm}_{\Phi}(\al)$. Then $\pi \in \mc{O}_L$ and $\sigma(\pi)\cdot \ov{\sigma(\pi)} = p$ for any $\sigma \in \absq$. Furthermore, the ideals of $\mc{O}_L$ generated by $\pi$ and $\ov{\pi}$ are coprime and hence, the abelian variety corresponding to the Weil $p$-integer $\pi$ is ordinary. 

Now, we have $\sigma_1(\al)\equiv r \Mod{\mfl}$, $\ov{\sigma}_2(\al)\equiv r^{1-m} \Mod{\mfl}$ and $\sigma(\al)\equiv 1\Mod{\mfl}$ for every $\sigma\in \absq\setminus \{\sigma_1^{-1},\ov{\sigma}_2^{-1} \}$. In particular, since $\sigma(\al)\equiv 1\Mod{\mfl}$ for every $\sigma\in \tau\Phi$, we have $\tau(\pi)\equiv 1\Mod{\mfl}$. Thus, $\ell$ divides the order $|B_{\pi}(\bFp)|$.

Furthermore, $p = \pi\cdot \ov{\pi} \equiv r\cdot r^{m-1}\equiv r^{m}\Mod{\mfl}$. So $p$ has order $n$ modulo $\mfl$ and hence, $B_{\pi}$ has embedding degree $n$ with respect to $\ell$. 

Lastly, $\pi\equiv r\Mod{\mfl}$ and hence, $\pi$ has order $m\cdot n$ modulo $\mfl$. Furthermore, any $\absq$-conjugate of $\pi$ is a product over some Galois conjugates of $\al$ and we have $\sigma(\al)^{m\cdot n}\equiv 1\Mod{\mfl}$ for every $\sigma\in \Gal(\wti{L}/\bq)$. Hence, $\sigma(\pi)^{m\cdot n}\equiv 1 \Mod{\mfl}$. Thus, $B_{\pi}$ has full embedding degree $m\cdot n$ with respect to $\ell$. \end{prf}

To construct such a Weil $p$-integer $\pi$, the most straightforward way would be as follows: 

\noindent 1. Randomly sample an algebraic integer $\al\in \mc{O}_{\wti{L}}$ that satisfy the congruences described in the proof.

\noindent 2. Subject $p:= \mr{Nm}_{\wti{L}/\bq}(\al)$ to a probabilistic primality test such as Miller-Rabin.

\noindent 3. If $p$ passes the primality test, accept $\pi:= \mr{Nm}_{\Phi}(\al)$ as a Weil $p$-integer with the desired properties. If not, return to Step 1 (i.e. sample another value for $\al$). 

The prime number theorem implies that the expected runtime for finding a suitable Weil number $\pi$ in this manner is polynomial in $\log(\ell)$.

\subsection{\fontsize{11}{11}\selectfont Jacobians with prescribed embedding and full embedding degrees}

\noindent We discuss the existence of Jacobians of small dimensions that fulfill the hypothesis of the preceding theorem. We will need the following two results regarding ordinary abelian surfaces from [How95].

\begin{Lem} $(\mr{[How95],\; Lemma\;} 12.1)$ Let $B_{\pi}$ be a simple ordinary abelian surface over a finite field $\bFq$ such that no abelian variety in its isogeny class over $\bFq$ is principally polarized. Then $\bq(\W_{B_{\pi}})$ is biquadratic. \end{Lem}

\begin{Thm} $(\mr{[How95],\; Theorem\;} 1.3)$ Any principally polarized simple ordinary abelian surface over a finite field $\bFq$ is a Jacobian of a hyperelliptic curve over $\bFq$.\end{Thm}

\vspace{0.5mm}

\begin{Corr} Let $L$ be a fixed CM field of degree $2g$ with $g\in \{2,3\}$, $m,n$ prescribed integers and $\ell$ a rational prime $\equiv 1 \pmod{m\cdot n}$ that splits completely in $L/\bq$. There exists a smooth projective curve $C$ of genus $g$ over a finite field $\bFq$ such that:

\noindent - $\End^0(\mr{Jac}(C))\cong L$.\\
- $\ell$ divides the order $\big|\mr{Jac}(C)(\bFq)\big|$ of the group of $\bFq$-points.\\
- $\mr{Jac}(C)$ has embedding degree $n$ with respect to $\ell$.\\
- $\mr{Jac}(C)$ has full embedding degree $n$ with respect to $\ell$, i.e. $[\bFq(\mr{Jac}(C)[\ell]):\bFq] = m\cdot n$.

\noindent If $g = 2$ and $L$ is not a biquadratic field, there exists a prime $p$ and a hyperelliptic curve $C$ over $\bFp$ fulfilling these conditions.\end{Corr}

\vspace{0.5mm}

\begin{prf} We construct a simple ordinary abelian variety $B$ of dimension $g$ over a prime field $\bFp$ fulfilling the conditions in the preceding theorem. We treat the cases $g=2$ and $g= 3$ separately.

\noindent \underline{Case 1}: $g= 2$ 

It follows from Theorem 13.3 of [How95] that there exists a Jacobian of a curve in the same isogeny class as $B$ over $\bFp$. Furthermore, this curve is hyperelliptic since it is smooth projective of genus two.

In particular, if $L$ is not biquadratic, Lemma 12.1 of [How95] implies that there exists a principally polarized abelian variety $B_1$ over $\bFp$ that is isogenous to $B$. Now, by ([How95], Theorem 1.3), $B_1$ is the Jacobian of a hyperelliptic curve $C_1$ over $\bFp$, which completes the proof.

\noindent \underline{Case 2}: $g= 3$ 

Since $B$ is simple of odd dimension, Theorem 1.2 of [How95] implies that there exists a principally polarized abelian $\wti{B}$ over $\bFp$ that is isogenous to $B$. It is well-known that every simple principally polarized abelian variety of dimension $\leq 3$ over an algebraically closed field is the Jacobian of some smooth projective curve. Hence, there exists a prime power $q$ such that $\wti{B}\times_{\bFp} \bFq$ is the Jacobian of a smooth projective genus three curve that is either hyperelliptic or planar quartic.\end{prf}

For a principally polarized abelian surface $B_{\pi}$ over a prime finite field $\bFp$, we have $\# B_{\pi} = P_{B_{\pi}}(1)$. On the other hand, if $C$ is a genus two hyperelliptic curve with $\Jac(C) = B$, then \vspace{-0.1cm}$$\# B_{\pi}(\bFp) = \frac{(\# C(\bFp))^2 + \# C(\bF _{p^2})}{2} - p.$$ Hence, to construct such a curve, it suffices to find a curve whose sets of $\bFp$ and $\bF _{p^2}$-points satisfy this equation. Since $\# C(\bFp))$ and $\# C(\bF _{p^2})$ lie in the respective Hasse-Weil intervals, we have $\# C(\bFp) = p+1-a_1$,  $\# C(\bF _{p^2}) = p^2+1+2a_2-a_1^2$ where $P_{B_{\pi}}(X) = X^4- a_1X^3 + a_2X^2 - pa_3 X+ p^2$. Algorithms for the construction of such curves using Gundlach invariants are studied in [LY10].

\section{\fontsize{11}{11}\selectfont Abelian varieties of type $\IV(1,d)$}


In this section, we study abelian varieties of type $\IV(1,d)$ over finite fields. The purpose here is to limit the size of the number field $\bq(\W_{B_{\pi}})$ while allowing the dimension of the (absolutely simple) abelian variety and the endomorphism algebra to be arbitrarily large. \vspace{1mm}

\begin{Def} \normalfont A simple abelian variety $A$ over a field $F$ is said to be of type $\IV$ in Albert's classification if it satisfies the following conditions:

\noindent 1. $\End^0(A)$ is a central division algebra over a CM field $K$ equipped with an involution $\;*\;$ of the second kind.\\
\noindent 2. $\End^0(A)$ is split at all the places $v$ of $K$ such that $v = {v}^*$.\\
\noindent 3. For any place $v$ of $K$, we have $\mr{inv}_v(\End^0(A))+\mr{inv}_{{v}^*}(\End^0(A)) = 0$ in $\bq/\bz$.\end{Def}

Furthemore, we say such an abelian variety $A$ is of type $\IV(e,d)$ if the center $K$ of $\End^0(A)$ is a CM field of degree $2e$ and $\End^0(A)$ is of dimension $d^2$ over its center. When the field of definition is a finite field, Honda-Tate theory implies that the dimension of $A$ is given by \vs $$ \dim A = \frac{1}{2}\cdot [K:\bq] \cdot [\End^0(A):K]^{1/2}  = e\cdot d. \vs $$

We say  $A$ is \textit{potentially} of type $\IV(e,d)$ if the base change to a suitable field extension $\wti{F}/F$ yields an abelian variety $\wti{A}:= A\times_F \wti{F}$ that is of type $\IV(e,d)$. In the paragraph following Corollary \ref{fieldOfDef}, we discuss a class of absolutely simple abelian varieties that are not of type $\IV(1,d)$ but are \textit{potentially} of type $\IV(e,d)$. In other words, a suitable base change yields an abelian variety with more endomorphisms.

Over a finite field of size $q$, a simple abelian $B_{\pi}$ corresponding to a Weil $q$-integer $\pi$ is of type $\IV(1,d)$ if and only if $\bq(\pi)$ is an imaginary quadratic field and the least common denominator of the Newton slopes of $B_{\pi}$ is $d$. We will use these abelian varieties to demonstrate that there exist absolutely simple abelian varieties of a prescribed dimension and embedding degree with respect to some rational prime dividing the size of the group of $\bFq$-points.

\vspace{1.5mm}

\begin{Lem} Let $g$ be an integer $\geq 3$. Let $B$ be a simple $g$-dimensional abelian variety over a finite field $\bFq$ such that one of the Newton slopes of $B$ is $\frac{j}{g}$ for some integer $j<g$ such that $\gcd(j, g)=1$. Then $B$ is absolutely simple.\end{Lem}

\begin{prf} Let $\pi$ be the Weil number corresponding to $B$. Suppose, by way of contradiction, that $\wti{B}= B\times_{\bFq} \bF _{q^N}$ is not simple. Since $\wti{B}$ has the same Newton polygon as $B$, it follows that $\wti{B}$ has a simple abelian subvariety $B'$ which has $\frac{j}{g}$ as a Newton slope. Let $m$ be the multiplicity of the slope $\frac{j}{g}$ in the Newton polygon of $B'$. Then $m < \dim B' < g$. Since a Newton polygon has integral breakpoints, $\frac{m\cdot j}{g}\in \bz$ and hence, $g$ divides $m$. Since $m < g$ we have a contradiction.\end{prf}

\vspace{1.5mm}

\begin{Prop} Let $B$ be a simple abelian variety of type $\IV(1, g)$ over a finite field $\bFq$ for some integer $g\geq 3$. Then $B$ has dimension $g$ and has Newton polygon $g\times \frac{j}{g}\;,\; g\times \frac{g-j}{g}$ for some integer $j< \frac{g}{2}$ prime to $g$.\end{Prop}

\begin{prf} It follows from the definition of an abelian variety of type $\IV(1,g)$ that the endomorphism algebra $\End^0(B)$ is a central division algebra dimension $g^2$ over its center $K$ which is an imaginary quadratic field. Thus, the dimension of $B$ is given by \vs $$\dim(B) = \frac{1}{2}\cdot [K:\bq]\cdot [D:K]^{1/2  } = g.$$  

Let $\pi$ be the Weil number associated to $B$. Then $\bq(\pi)$ is an imaginary quadratic field in which $p$ splits. Let $v$ be any $p$-adic valuation on $\bq(\pi)$. The Newton slopes of $B$ are given by \vspace{-1mm}$$g\times \frac{v(\pi)}{v(q)}\;\;,\;\; g\times \frac{v(q)-v(\pi)}{v(q)}\vspace{-1mm} $$ and since $p$ splits in $\bq(\pi)/\bq$, the set of Hasse invariants of $\End^0(B)$ is the same as the set of Newton slopes of $B$. Now, $g = [\End^0(B):\bq(\pi)]$ is the least common denominator of the Hasse invariants. Hence, the least common denominator of the Newton slopes is $g$, meaning that $\frac{v(\pi)}{v(q)} = \frac{j}{g}$ for some integer $j< \frac{g}{2}$ prime to $g$.\end{prf}

\vspace{1.5mm}

\begin{Corr} Any abelian variety of type $\IV(1,d)$ over a finite field is absolutely simple.\end{Corr}
\begin{prf} This follows immediately from the last two propositions.\end{prf}

\vspace{1mm}

\begin{Prop} Let $B_{\pi}$ be an absolutely simple abelian variety over a finite field $\bFq$ such that $|\W_{B_{\pi}}|\leq 2$. Then $B$ is either an elliptic curve or an abelian variety of type $\IV(1,d)$ where $d = \dim B$.\end{Prop}

\begin{prf} Any supersingular abelian variety over the algebraic closure $\ov{\bF} _q$ is isogenous to a a power of a supersingular elliptic curve over $\ov{\bF} _q$. Thus, if $B_{\pi}$ is not an elliptic curve, it can be assumed to not be supersingular.

Clearly, $[\bq(\pi):\bq] = \big|\W_{B_{\pi}}\big|\leq 2$. The only real Weil $q$-integers are $\pm \sqrt{q}$. So, if $\bq(\pi)$ is a totally real field, then $B_{\pi}$ is supersingular. Hence, we may assume $\bq(\pi)$ is a imaginary quadratic field. Write $D_{\pi}:=\End^0(B_{\pi})$. Then \vspace{-0.15cm} $$d = \dim B_{\pi} = [D_{\pi}:\bq(\pi)]^{1/2}$$ and hence, $[D_{\pi}:\bq(\pi)] = d^2$.\end{prf}

\vspace{1.5mm}

\begin{Prop} Let $B_{\pi}$ be an abelian variety of type $\IV(1,d)$ over $\bFq$ and $\ell$ be a prime such that $\ell$ divides $\big|B_{\pi}(\bFq)\big|$ and $\ell\nmid q\cdot (q-1)$. Then the embedding degree of $B_{\pi}$ with respect to $\ell$ coincides with the full embedding degree $[\bFq(B_{\pi}[\ell]):\bFq]$.\end{Prop}

\begin{prf} Since $\ell$ divides $\big|B_{\pi}(\bFq)\big|$, there exists a prime $\mf{l}$ of $\bq(\pi)$ lying over $\ell$ such that $\pi\equiv 1\pmod{\mf{l}}$. Let $N$ be the embedding degree. Then $q$ has order $N$ in $(\mc{O}_{\bq(\pi)}/\mf{l})^*$ and since $q\equiv \ov{\pi} \pmod{\mf{l}}$, it follows that $\ov{\pi}$ has order $N$ in $(\mc{O}_{\bq(\pi)}/\mf{l})^*$. Thus, $\pi^N\equiv 1\pmod{\ell \mc{O}_{\bq(\pi)}}$, which implies that the degree $[\bFq(B_{\pi}[\ell]):\bFq]$ divides $N$.\end{prf}

\vspace{1.5mm}

\begin{Prop} Let $p$ be a fixed prime and let $d\geq 3$, $N$ $\geq 2$  be prescribed integers. Then there exists an absolutely simple abelian variety $B$ of dimension $d$ over a finite field $\bFq$ of characteristic $p$ and a prime divisor $\ell$ of $|B(\bFq)|$ such that 

\noindent - $B$ has embedding degree $N$ with respect to $\ell$.

\noindent - For any extension $F$ of $\bFq$ and any prime $\ell'$ dividing $|B(F)|$, the embedding degree and the full embedding degree of $B$ with respect to $\ell'$ coincide.\end{Prop}

\begin{prf} Let $B$ be an absolutely simple abelian variety of type $\IV(1,d)$ over $\bFq$. Then $\dim B = d$ and $L:=\bq(\pi)$ is an imaginary quadratic field in which $p$ splits. By construction, the subgroup of $L^*/L^{*N}$ generated by $\{\pi,\ov{\pi}\}$ is isomorphic to $(\bz/N\bz)^2$. From Kummer theory, it follows that $$\Gal(L(\zeta_{_N}, \sqrt[N]{\pi},\sqrt[N]{\ov{\pi}})/L(\zeta_{_N}))\cong (\bz/N\bz)^2.$$

Choose a prime $\mf{l}$ of $L(\zeta_{_N})$ such that:

\noindent - $\mf{l}$ has local degree one over $\bq$.

\noindent - $\mf{l}$ splits in $L(\zeta_{_N}, \sqrt[N]{\pi})$. 

\noindent - $\mf{l}$ is inert in $L(\zeta_{_N}, \sqrt[N]{\ov{\pi}})$. 

The set of primes fulfilling these conditions has Dirichlet density ${\phi(N)}\cdot {N^{-2}}$ and in particular, such a prime exists. Let $\ell$ be the rational prime lying under $\mf{l}$. Set $\pi_{_0}:= \pi^{^{\frac{\ell-1}{N}}}$, $q_{_0}:= q^{^{\frac{\ell-1}{N}}}$. Then $\pi_{_0}\equiv 1 \pmod{\mf{l}}$ and hence, $\ell$ divides $|B(\bF _{q_{_0}})|$. On the other hand, $q_{_0}$ has order $N$ in $(\mc{O}_{L(\zeta_{_N})}/\mf{l})^*$ and hence, $\ell$ divides $\Phi_N(q_{_0})$. Thus, the abelian variety $B_{\pi_{_0}} = B_{\pi}\times_{\bFq} \bF _{q_{_0}}$ and the prime $\ell$ fulfill the conditions. \end{prf}

\noindent \textbf{Example.} We provide an example of an abelian variety of type $\IV(1,d)$ that arises as a Jacobian of a hyperelliptic curve. Let $d$ be an odd prime such that $\ell := 2d+1 $ is also a prime. Consider the genus $\ell$ hyperelliptic curve \vspace{-0.15cm} $$C:Y^2 = 1- X^{\ell}$$ over $\bq$. Let $\zeta_{\ell}\in \qbar$ be a fixed primitive $\ell$-th root of unity. The curve $C\times_{\bq}{\bq(\zeta_{\ell})}$ is endowed with the automorphism $$(x_1, y_1)\mapsto (\zeta_{\ell}\cdot x_1,\; y_1) $$ and hence, the Jacobian $\Jac(C)$ has action by the ring $\bz[\zeta_{\ell}]$. Since \vs $$[\bq(\zeta_{\ell}):\bq] = \ell-1 = \dim \Jac(C),$$ it follows that $\Jac(C)$ is simple with endomorphism ring $\bz[\zeta_{\ell}]$. 

The curve $C$ has good reduction away from the set $\{2, \ell \}$. Since the extension $\bq(\zeta_{\ell})/\bq$ is abelian, it follows from the theory of complex multiplication that the Newton slopes of the reduction $\Jac(C)_p$ at any prime of good reduction are determined by the splitting of $p$ in the field $\bq(\zeta_{\ell})$. In particular, if $p$ is a prime of inertia degree $d$ in $\bq(\zeta_{\ell})$, the Newton polygon of $\Jac(C_p)$ is given by \vspace{-1mm} $$d\times \frac{j}{d}\;\; ,\;\; d\times \frac{d-j}{d}  $$ for some integer $j\leq d-1$ prime to $d$. Thus, $\Jac(C_p)$ is absolutely simple and $\Jac(C_p)\times_{\bFp}\bF _{p^d}$ is of type $\IV(1,d)$.

\subsubsection{\fontsize{11}{11}\selectfont The minimum field of definition}

\begin{Lem} Let $B_{\pi}$ be an absolutely simple abelian variety over a finite field $k$ with\\ $[D_{\pi}:\bq(\pi)] = d^2$. Suppose the extension $\bq(\W_{B_{\pi}})/\bq$ is such that the decomposition group of a prime of $\bq(\W_{B_{\pi}})$ lying over $p$ is a normal subgroup of $\Gal(\bq(\W_{B_{\pi}})/\bq)$. Then:

\noindent $1$. $p$ splits completely in $\bq(\W_{B_{\pi}})/\bq$.\\
$2$. The Newton slopes of $B_{\pi}$ have least common denominator $d$.\end{Lem}

\begin{prf} This is lemma 5.3 of [Th17].\end{prf}

\noindent \textbf{Remark:} In particular, the hypothesis of the lemma is satisfied when $\bq(\pi)/\bq$ is abelian.

\begin{Corr} \label{fieldOfDef} If there exists an abelian variety $B$ of type $\IV(1,d)$ over $\bF _{p^k}$, then $d$ divides $k$.\end{Corr}

\begin{prf} Let $\pi$ be the corresponding Weil number. As noted earlier, $B$ is \textit{absolutely} simple. Furthermore, the extension $\bq(\pi)/\bq$ is quadratic and in particular, abelian. Hence, it follows from the lemma that $p$ splits (completely) in $\bq(\pi)/\bq$. From Proposition 2.14, it follows that the Newton polygon is \vspace{-1mm}$$d\times\frac{j}{d}\;\;,\;\; d\times\frac{d-j}{d} \vspace{-1mm} $$ for some integer $j<d$ such that $\gcd(j,d) = 1$. The Newton slopes are given by $\frac{v(\pi)}{v(p^k)} = \frac{v(\pi)}{k\cdot v(p)}$ where $v$ runs through the places of $\bq(\pi)$ lying over $p$. Since $p$ splits completely in $\bq(\pi)/\bq$, it follows that $\frac{v(\pi)}{v(p)}$ is an integer and hence, have $d$ divides $k$.\end{prf}

Note that the abelian variety $B$ could have a $B_0$ model over a proper subfield of $\bF _{p^d}$, even though $B_0$ would not be of type $\IV(1,d)$.\vspace{1mm} 

\noindent \textbf{Example.} For instance, let $\pi_1$ be an ordinary Weil $p$-integer corresponding to an ordinary elliptic curve over $\bFp$ and set $\pi:= \sqrt[d]{\pi_1^{d-1}\cdot\ov{\pi}_1}$. Then $\pi$ is a Weil $p$-integer corresponding to a simple abelian variety of dimension $d$ over $\bFp$ and with the endomorphism algebra $\bq(\pi)$ a CM field of degree $2d$. On the other hand, the Weil $p^d$-integer $\pi^d = \pi_1^{d-1}\cdot\ov{\pi}_1$ is an algebraic integer of degree $2$ and the corresponding abelian variety over $\bF _{p^d}$ is of type $\IV(1,d)$.\\

\noindent We now show that the converse to Corollary \ref{fieldOfDef} is also true.

\begin{Prop} Let $p$ be a prime and $d$ any integer $\geq 3$. There exists an abelian variety of type $\IV(1,d)$ over $\bF _{p^d}$.\end{Prop}

\begin{prf} Let $\pi$ be a Weil $p$-integer corresponding to an ordinary elliptic curve over $\bFp$, meaning that $\gcd(\pi,\ov{\pi}) = (1)$ and $p = \pi\cdot\ov{\pi}$. Then $\pi$ and $\ov{\pi}$ generate distinct principal prime ideals in $\bq(\pi)$. 

Now, $\wti{\pi}:= \pi^{d-1}\cdot\ov{\pi}$ is a Weil $p^d$-integer. Let $B_{\wti{\pi}}$ be the corresponding abelian variety (up to isogeny) over $\bF _{p^d}$. The Newton polygon of $B_{\wti{\pi}}$ is \vs $$d\times\frac{1}{d}\;\;,\;\; d\times\frac{d-1}{d}.$$ Clearly, $\bq(\wti{\pi}) = \bq(\pi)$ and in particular, the center of $\End^0(B_{\wti{\pi}})$ is the imaginary quadratic field $\bq(\pi)$. Furthermore, there are two primes in $\bq(\pi)$ lying over $p$ and the Hasse invariants of $\End^0(B)$ are given by $$\mr{inv}_v(\End^0(B_{\wti{\pi}}))=\begin{cases}
0 & \text{ if } v\nmid p\\
\frac{v(\wti{\pi})}{v(q)}\cdot  [\bq(\pi)_v:\bq_p]\;\pmod{1} & \text{ if } v|p.\end{cases}$$

Since $p$ splits (completely) in $\bq(\pi)/\bq$, we have $[\bq(\pi)_{(\wti{\pi})}:\bq_p] = [\bq(\pi)_{(\ov{\wti{\pi}})}:\bq_p] = 1$ and hence, 
$\End^0(B_{\wti{\pi}})$ is ramified at the primes $(\pi)$ and $(\ov{\pi})$ with Hasse invariants $\frac{1}{d}$ and $\frac{d-1}{d}$ respectively. In particular, the least common denominator of the Hasse invariants is $d$ and hence, $\End^0(B_{\wti{\pi}})$ is of dimension $d^2$ over its center.\end{prf}

\bigskip

\noindent \textbf{Acknowledgements:} The author thanks the anonymous referee for helpful feedback and for suggesting several improvements to previous drafts.

\bigskip 

\begin{center}\textbf{References} \end{center}
\small

\noindent [BCF] Benger, Charlemagne, D. Freeman, \textit{On the security of pairing-friendly abelian varieties over non-prime fields} \vspace{0.1cm}

\noindent [FFS08] D. Freeman, M. Streng, P. Stevenhagen, \textit{Abelian varieties with prescribed embedding degree}, Algorithmic Number Theory. ANTS 2008. Lecture Notes in Computer Science, vol 5011 \vspace{0.1cm}

\noindent [Gal] S. Galbraith, \textit{Supersingular curves in Cryptography} \vspace{0.1cm}

\noindent [GMV07] S. Galbraith, Mckee, Valenca, \textit{Ordinary abelian varieties having small embedding degree}, Finite Fields and their Applications, 13:800–814, 2007 \vspace{0.1cm}

\noindent [Hit] L. Hitt, \textit{On an Improved Definition of the Embedding Degree} \vspace{0.1cm}

\noindent [How95] E. Howe, \textit{Principally polarized ordinary abelian varieties over finite fields}, Transactions of the AMS, Volume 347, Number 7, July 1995 \vspace{0.1cm}

\noindent [Kow03] E. Kowalski, \textit{Some local-global applications of Kummer theory}, Manuscripta Math., 111(1), 2003.

\noindent [Kow06] \underline{\;\;\;\;\;\;\;}, \textit{Weil numbers generated by other Weil numbers}, J. London Math. Soc. (2) 74 (2006). \vspace{0.1cm}

\noindent [LY10] K. Lauter, T. Yang, \textit{Computing genus 2 curves from invariants on the Hilbert moduli space}  \vspace{0.1cm}

\noindent [MOV93] A. Menezes, T. Okamoto, and S. Vanstone. \textit{Reducing elliptic curve logarithms to logarithms in a finite field} IEEE Transactions on Information Theory 39 (1993),1639–1646. \vspace{0.1cm}

\noindent [Oort] F. Oort, \textit{Abelian Varieties over finite fields} \vspace{0.1cm}

\noindent [Shi60] G. Shimura, \textit{Abelian Varieties with Complex Multiplication and Modular Functions}, Princeton University Press \vspace{0.1cm}

\noindent [ST66] J.P. Serre, J. Tate, \textit{Good reduction of abelian varieties}, Ann. of Math. (2), 88, 1968 \vspace{0.1cm}

\noindent [Tat69] J. Tate, \textit{Classes d’isogenie des varietes abeliennes sur un corps fini}\vspace{0.1cm}

\noindent [Th17] S. Thakur, \textit{On some abelian varieties of type $IV$}, Preprint (in revision at Math Zeit.)

\bigskip

\normalsize
\noindent Email: stevethakur01@gmail.com

\noindent Steve Thakur\\
Cybersecurity Research Group\\
Battelle Institute\\
Columbus, OH.

\newpage

\bigskip
\bigskip

\end{document}